\documentclass[a4paper,12pt]{amsart}
\usepackage{amsfonts}
\usepackage{amssymb}
\usepackage{ifthen}
\usepackage{graphicx}
\nonstopmode \numberwithin{equation}{section}
\usepackage[usenames]{color}
\usepackage{enumerate}
\newtheorem{thm}{Theorem}[section]
\newtheorem{lem}[thm]{Lemma}
\newtheorem{cor}[thm]{Corollary}
\newtheorem{prop}[thm]{Proposition}

\newtheorem{cl}{Claim}[section]
\newtheorem{ca}{Case}[section]
\newtheorem{sca}{Subcase}[section]
\newtheorem{scl}[section]{Subclaim}
\newtheorem{conj}[equation]{Conjecture}

\theoremstyle{definition}
\newtheorem{defn}[thm]{Definition}
\newtheorem{op}[thm]{Open Problem}
\newtheorem{ques}[thm]{Question}
\newtheorem{rem}[thm]{Remark}
\newtheorem{exam}[equation]{Example}

\newcounter {own}
\def\theown {\thesection       .\arabic{own}}

\newenvironment{pf}[1][]{%
 \vskip 3mm
 \noindent
 \ifthenelse{\equal{#1}{}}%
  {{\slshape Proof. }}%
  {{\slshape #1.} }%
 }%
{\qed\bigskip}

\newcounter{alphabet}

\newcommand{\IR}{{\mathbb R}}

\newcommand{\diam}{{\operatorname{diam}}}

\newcommand{\dist}{{\operatorname{dist}}}

\def\be{\begin{equation}}
\def\ee{\end{equation}}

\newcommand{\ben}{\begin{enumerate}}
\newcommand{\een}{\end{enumerate}}

\newcommand{\blem}{\begin{lem}}
\newcommand{\elem}{\end{lem}}
\newcommand{\bthm}{\begin{thm}}
\newcommand{\ethm}{\end{thm}}
\newcommand{\bcor}{\begin{cor}}
\newcommand{\ecor}{\end{cor}}
\newcommand{\beg}{\begin{exam}}
\newcommand{\eeg}{\end{exam}}
\newcommand{\begs}{\begin{examples}}
\newcommand{\eegs}{\end{examples}}
\newcommand{\bdefe}{\begin{defn}}
\newcommand{\edefe}{\end{defn}}
\newcommand{\bprob}{\begin{prob}}
\newcommand{\eprob}{\end{prob}}
\newcommand{\bques}{\begin{ques}}
\newcommand{\eques}{\end{ques}}
\newcommand{\bei}{\begin{itemize}}
\newcommand{\eei}{\end{itemize}}
\newcommand{\bcon}{\begin{conj}}
\newcommand{\econ}{\end{conj}}
\newcommand{\bop}{\begin{op}}
\newcommand{\eop}{\end{op}}

\newcommand{\bca}{\begin{ca}}
\newcommand{\eca}{\end{ca}}
\newcommand{\bsca}{\begin{sca}}
\newcommand{\esca}{\end{sca}}

\newcommand{\bcl}{\begin{cl}}
\newcommand{\ecl}{\end{cl}}

\newcommand{\bscl}{\begin{scl}}
\newcommand{\escl}{\end{scl}}

\newcommand{\bcons}{\begin{conjs}}
\newcommand{\econs}{\end{conjs}}
\newcommand{\bprop}{\begin{propo}}
\newcommand{\eprop}{\end{propo}}
\newcommand{\br}{\begin{rem}}
\newcommand{\er}{\end{rem}}
\newcommand{\brs}{\begin{rems}}
\newcommand{\ers}{\end{rems}}
\newcommand{\bo}{\begin{obser}}
\newcommand{\eo}{\end{obser}}
\newcommand{\bos}{\begin{obsers}}
\newcommand{\eos}{\end{obsers}}
\newcommand{\bpf}{\begin{pf}}
\newcommand{\epf}{\end{pf}}
\newcommand{\ba}{\begin{array}}
\newcommand{\ea}{\end{array}}
\newcommand{\beq}{\begin{eqnarray}}
\newcommand{\beqq}{\begin{eqnarray*}}
\newcommand{\eeq}{\end{eqnarray}}
\newcommand{\eeqq}{\end{eqnarray*}}

\newcounter{minutes}
\divide\time by 60
\newcounter{hours}
\multiply\time by 60 \addtocounter{minutes}{-\time}

\begin{document}

\bibliographystyle{amsplain}
\title{Gromov hyperbolic John is quasihyperbolic John II}

\author{Qingshan Zhou${}^{\mathbf{*}}$}
\address{Qingshan Zhou, School of Mathematics and Big Data, Foshan university,  Foshan, Guangdong 528000, People's Republic
of China} \email{qszhou1989@163.com; q476308142@qq.com}

\author{Saminathan Ponnusamy}
\address{Saminathan Ponnusamy, Department of Mathematics, Indian Institute of Technology Madras, Chennai 600036,
India}
\email{samy@iitm.ac.in}

\def\thefootnote{}
\footnotetext{ \texttt{\tiny File:~\jobname .tex,
          printed: \number\year-\number\month-\number\day,
          \thehours.\ifnum\theminutes<10{0}\fi\theminutes}
} \makeatletter\def\thefootnote{\@arabic\c@footnote}\makeatother

\date{}
\subjclass[2000]{Primary: 30C65, 30L10; Secondary: 30C20  30F45} \keywords{John spaces, quasihyperbolic geodesic, Gromov hyperbolic spaces.\\
${}^{\mathbf{*}}$ Corresponding author}

\begin{abstract} In this paper, we introduce a concept of quasihyperbolic John spaces (with center) and provide a criteria to determine spaces to be quasihyperbolic John. As an application, we provide a simple proof to show that a John space with a Gromov hyperbolic quasihyperbolization is quasihyperbolic John, quantitatively. This gives an affirmative answer to an open question posed by Heinonen (Rev.~Math.~Iber, 1989), which has been studied by Gehring et al. (Math.~Scand, 1989).
\end{abstract}

\thanks{The research was partly supported by NSF of
China (No. 11901090).}

\maketitle{} \pagestyle{myheadings} \markboth{Gromov hyperbolic John is quasihyperbolic John II}{Zhou Qingshan and Saminathan Ponnusamy}

\section{Introduction and main results}

Throughout this paper, the set $X$ is a locally compact, incomplete and rectifiably connected metric space $(X, d)$, and $x_0\in X$ is a base point. Every proper domain of Euclidean spaces $\mathbb{R}^n$ is such a space. We focus our attention on quasihyperbolic geometry of John spaces. The class of John domains in $\mathbb{R}^n$ was first considered by John \cite{Jo61} in the study of elasticity theory. For $a\geq 1$, the space $X$ is called $a$-{\it John with center $x_0\in X$} if for all $x\in X$, there is an arc $\gamma$ in $X$ joining $x$ to $x_0$ such that for all $z\in \gamma$,
$$\ell(\gamma[x, z])\leq a\,d(z),
$$
where $\gamma[x, z]$ is the subarc of $\gamma$ with endpoints $x$ and $z$, and
$d(z)$ denotes the distance between $z$ and $\partial X$,  the metric boundary of $X$.
The arc $\gamma$ is called an {\it $a$-cone arc}. Moreover, $X$ is called {\it John} if there is a constant $a\geq 1$  and a distinguished point $x_0\in X$ such that $X$ is $a$-John with center $x_0\in X$. Note that John spaces are bounded.

The quasihyperbolic metric $k$ in $X$ is defined by
$$k(x, y)=\inf \Big\{\int_{\gamma} \frac{|dz|}{d(z)}\Big\},
$$
where the infimum is taken over all rectifiable curves $\gamma$ in $X$ with the endpoints $x$ and $y$, and $|dz|$ denotes the arc-length element with respect to the metric $d$. We remark that the resulting space $(X,k)$ is complete, proper and geodesic (cf. \cite{BHK}).
This notation was introduced by Gehring and Palka \cite{GP76} in Euclidean setting, and for further geometric properties, see \cite{GO}. The importance of the quasihyperbolic metric and John domains in quasiconformal mapping theory is well understood, see for example \cite{BKL,CP17,GNV94,Guo15,HLPW15,NV}.

In \cite{GHM}, Gehring et al. studied geometric properties of quasihyperbolic geodesics in planar domains and demonstrated that quasihyperbolic geodesics in simply connected planar John domains are cone arcs. This result is not valid in multi-connected domains, for examples see \cite{GHM}. Naturally, it raises the problem of determining sufficient and/or necessary conditions in a John domain such that quasihyperbolic geodesics to be cone arcs. This problem has been considered by several authors. See \cite{BHK,GHM,Guo15,Hei89,ZLR}  for further details and background information. In particular, Heinonen asked in \cite[Question 2]{Hei89} whether quasihyperbolic geodesics are cone arcs in John domains of $\IR^n$ which quasiconformally equivalent to the unit ball or uniform domains.

Note that the previous results in $\mathbb{R}^n$ depend on $n$. Also, the proofs are highly dependent on Euclidean geometry. In this paper, we investigate this problem by introducing the following concept in metric spaces setting and search for dimension-free results.

\begin{defn} Let $x_0\in X$ and $a\geq 1$. The space $X$ is called {\it quasihyperbolic $a$-John with center $x_0$} if for each $x\in X$, every quasihyperbolic geodesic $\alpha$ in $X$ joining $x$ to $x_0$ is an $a$-cone arc. Moreover, $X$ is called {\it quasihyperbolic John} if there is a constant  $a\geq 1$ and a point $x_0\in X$ such that $X$ is quasihyperbolic $a$-John with center $x_0$.
\end{defn}

Clearly, every quasihyperbolic John space is John. There are John spaces which are not quasihyperbolic John, see \cite[Examples 5]{GHM}. Indeed, we shall see that the class of quasihyperbolic John spaces (with center) is very wide.

As our first main result, the following theorem  gives a sufficient and necessary condition  in determining which space is quasihyperbolic John.

\begin{thm}\label{thm-0} Let $x_0\in X$. The space $X$ is quasihyperbolic $a$-John with center $x_0$ if and only if
\begin{enumerate}
  \item There is a constant $A\geq 1$ such that $A d(x_0)\geq \diam (X)$, and
  \item For every $x\in X$ and for each quasihyperbolic geodesic $\alpha$ between $x$ and $x_0$, the following holds: If $y_1$ and $y_2$ are points on $\alpha$ satisfying for each $u\in \alpha[y_1,y_2]$, $d(u)\leq 2\min\{d(y_1),d(y_2)\}$, then we have
      $$k(y_1,y_2)\leq A.$$
\end{enumerate}
The constants $a$ and $A$ depending only on each other.
\end{thm}

Using Theorem \ref{thm-0}, one easily finds that bounded uniform spaces are quasihyperbolic John. Note that uniform metric spaces were introduced by Bonk et al. in \cite{BHK}, where they proved quasihyperbolic geodesics in uniform spaces are uniform arcs. Recall that $X$ is called $c$-{\it uniform} for $c\geq 1$ if each pair of points $x$, $y$ in $X$ can
be joined by a rectifiable curve $\gamma$ in $X$ satisfying

$(a)$ $($Quasiconvexity condition$)$ $\ell(\gamma)\leq c\,d(x,y)$, and

$(b)$ $($Double cone condition$)$ $\min\{\ell(\gamma[x,z]),\ell(\gamma[z,y])\} \leq c\,d(z)$ for all $z\in \gamma$.

Moreover, $\gamma$ is called a $c$-uniform curve.

\begin{prop}  A bounded $c$-uniform space $X$ is quasihyperbolic $a$-John with $a$ depending only on $c$.
\end{prop}
\bpf We only need to verify the criteria given in Theorem \ref{thm-0}. For (1), since $X$ is bounded, we may choose a point $x_0\in X$ with $d(x_0)=\max\{d(z): z\in X\}$. For every $x\in X$, there is a $c$-uniform curve $\gamma$ connecting $x$ to $x_0$. Let $y_0\in \gamma$ be such that $\ell(\gamma[x,y_0])=\ell(\gamma[y_0,x_0])$. The uniformity of $\gamma$ leads to
$$d(x,x_0)\leq \ell(\gamma)\leq 2cd(y_0)\leq 2cd(x_0),$$
which yields $
\diam (X)\leq 4c d(x_0)$ by using the arbitrariness of $x$ in $X$. Hence (1) is true.

For (2), let $x\in X$ and let $\alpha$ be a quasihyperbolic geodesic  between $x$ and $x_0$. Suppose that $y_1$ and $y_2$ are points on $\alpha$ satisfying for each $u\in \alpha[y_1,y_2]$, $d(u)\leq 2\min\{d(y_1),d(y_2)\}$. By \cite[Theorem 2.10]{BHK}, we know that the quasihyperbolic geodesic $\alpha[y_1,y_2]$ is $c_1$-uniform with $c_1=c_1(c)$. Take a point $y\in \alpha[y_1,y_2]$ be such that $\ell(\alpha[y_1,y])=\ell(\alpha[y,y_2])\leq c_1d(y)$. Since $d(y)\leq 2\min\{d(y_1),d(y_2)\}$, we obtain
$$d(y_1,y_2)\leq \ell(\alpha[y_1,y_2])\leq 2c_1d(y)\leq c_1\min\{d(y_1),d(y_2)\}.$$
This inequality, together with \cite[(2.16)]{BHK}, ensures that $$k(y_1,y_2)\leq 4c^2\log\Big(1+\frac{d(y_1,y_2)}{\min\{d(y_1),d(y_2)\}}\Big)\leq 4c^2\log(1+c_1).$$
Therefore, $(2)$ holds.  Hence, Theorem \ref{thm-0} implies that $X$ is quasihyperbolic $a$-John with $a$ depending only on $c$.
\epf

As the second application, we use the criteria of Theorem \ref{thm-0} to prove that Gromov hyperbolic John spaces are quasihyperbolic John. Here and hereafter, $X$ is called {\it Gromov hyperbolic} if it admits a {\it Gromov hyperbolic quasihyperbolization}. That is to say, $(X,k)$ is $\delta$-hyperbolic for some constant $\delta\geq 0$, where $k$ is the quasihyperbolic metric of $X$. Our second main result is as follows.

\begin{thm}\label{thm-1} Let $(X,d)$ be an $a$-John space with center $x_0\in X$, and let $(X,k)$ be Gromov $\delta$-hyperbolic. Then every quasihyperbolic geodesic in $X$ starting from $x_0$ is an $M$-cone arc with $M=M(a,\delta)$.
\end{thm}

\br
\begin{enumerate}
  \item Theorem \ref{thm-1} asserts that John spaces with  center carrying a Gromov hyperbolic quasihyperbolization are quasihyperbolic John, quantitatively. So it answers affirmatively to \cite[Question 2]{Hei89} with the parameter $M$ independent of the dimension $n$.  Note that the class of Gromov hyperbolic John spaces includes uniform domains and inner uniform domains of $\IR^n$, finitely connected John domains in the plane, and John domains of $\mathbb{R}^n$ which are Gromov $\delta$-hyperbolic in the quasihyperbolic metric, etc. There are many applications of these domains in quasiconformal mapping theory and potential analysis, see \cite{CP17,GNV94,Guo15}.
  \item Theorem \ref{thm-1} is an improvement of \cite[Proposition 7.12]{BHK}, \cite[Theorem 4.1]{GHM} and  \cite[Remark 3.10]{Guo15}. These results were considered in bounded domains of Euclidean spaces $\mathbb{R}^n$, where the parameter $M$ depends also on $n$. Indeed, their proofs may not hold without assuming the spaces to be Ahlfors regular. This is because we do not know whether the spaces/domains satisfy the ball separation condition in this setting, see \cite{Guo15}.
  \item We do not need the space $X$ in Theorem \ref{thm-1} to be roughly starlike. So Theorem \ref{thm-1} gives an improvement of \cite[Theorem 1.1]{ZLR}. Moreover, our method of proof is more direct but simple than that of \cite[Theorem 1.1]{ZLR}. Another important observation is that every cone arc starting from $x_0$ is a quasihyperbolic quasigeodesic, see Lemma \ref{z2}. Using this observation and the stability property of Gromov hyperbolic spaces, one may easily check the criteria for quasihyperbolic John spaces as mentioned in Theorem \ref{thm-0}.
\end{enumerate}
\er

This paper is organized as follows. Section \ref{sec-2} contains notations and the basic definitions. The proofs of Theorems \ref{thm-0} and \ref{thm-1} are given in Section \ref{sec-3}.

\section{Preliminaries}\label{sec-2}
Let $(X, d)$ denote a metric space. Its metric completion is $\bar{X}$. The space $X$ is incomplete if its metric boundary  $\partial X$
is nonempty, where $\partial X=\bar{X}\setminus X$. For all $x\in X$, $d(x)=\dist(x,\partial X)$ if $\partial X\neq \emptyset$. The open ball with center $x\in X$ and radius $r>0$ is denoted by
$B(x,r)=\{z\in X:\; d(z,x)<r\}$.
A {\it curve} in $X$ is a continuous map $\gamma:\; I\to X$ from an interval $I\subset \IR$ to $X$. If $\gamma$ is an embedding of $I$, then it is also called an {\it arc}. We also denote the image set $\gamma(I)$ of $\gamma$ by $\gamma$. The {\it length} $\ell(\gamma)$ of $\gamma$ with respect to the metric $d$ is defined in an obvious way. Here the parameter interval $I$ is allowed to be open or half-open. A metric space $(X, d)$ is called {\it rectifiably connected} if every pair of points in $X$ can be joined with a
curve $\gamma$ in $X$ with $\ell_d(\gamma)<\infty$.

A {\it geodesic} arc $\gamma$ joining $x\in X$ to $y\in X$ is a map $\gamma:I=[0,l]\to X$ such that $\gamma(0)=x$, $\gamma(l)=y$ and
$d(\gamma(t),\gamma(t'))=|t-t'|$  for all $t,t'\in I.$
A metric space $X$ is said to be {\it geodesic} if every pair of points can be joined by a geodesic arc. We let $[x,y]$ be a geodesic between two points $x$ and $y$ in $X$. Let $\lambda\geq 1$ and $\mu\geq 0$.
A {\it $(\lambda, \mu)$-quasigeodesic} curve in $X$ is a $(\lambda,\mu)$-quasi-isometric embedding $\gamma:I\to X$, $I\subset \mathbb{R}$. More explicitly,
$${\lambda}^{-1}|t-t'|-\mu\leq d(\gamma(t),\gamma(t')) \leq \lambda|t-t'|+\mu ~\mbox{  for all $t,t'\in I.$}
$$

Let $\delta\geq 0$. A geodesic space $(X,d)$ is called (Gromov) $\delta$-hyperbolic if for all triples of geodesics $[x,y], [y,z], [z,x]$ in $(X,d)$, every point in $[x,y]$ is within the distance $\delta$ from $[y,z]\cup [z,x]$. Finally, we recall the stability property of quasigeodesics in Gromov hyperbolic spaces.

\begin{thm}\label{z1} $($\cite[Chapter III.H Theorem $1.7$]{BH99}$)$
For all $\delta\geq 0, \lambda\geq 1, \mu\geq 0$, there is a number $R=R(\delta,\lambda,\mu)$ with the following property: If $X$ is a $\delta$-hyperbolic space, $\gamma$ is a $(\lambda,\mu)$-quasigeodesic in $X$ and $[x,y]$ is a geodesic joining the endpoints of $\gamma$, then the Hausdorff distance between $[x,y]$ and $\gamma$ is less than $R$.
\end{thm}

\section{Proofs of Main results}\label{sec-3}
\subsection{A criteria for quasihyperbolic John spaces} Our goal of this subsection is to prove Theorem \ref{thm-0}. We start with some auxiliary results.
We first need  the following well-known inequalities (cf. \cite{BHK,ZLR}): for all $x,y\in X$
\be\label{zt-1} \Big|\log\frac{d(x)}{d(y)} \Big| \leq k(x,y)\ee
and
\be\label{zt-2} \log\Big(1+\frac{\ell(\gamma)}{\min\{d(x),d(y)\}}\Big) \leq \ell_k(\gamma),\ee
where $\gamma$ is a curve in $X$ with $x,y\in \gamma$ and $\ell_k(\gamma)$ is the quasihyperbolic length of $\gamma$.

\begin{lem}\label{t-1}
If $\gamma$ is an $a$-cone arc in $X$ joining $x$ to $x_0$, then for  each $u\in \gamma$ we have  $d(x)\leq 2a d(u).$
\end{lem}
\bpf Fix $u\in \gamma$. If $u\in B(x,d(x)/2)$, then we have
$$d(u)\geq d(x)-d(x,u)\geq d(x)/2.$$
Otherwise, by the assumption we obtain
$$ad(u)\geq \ell(\gamma[x,u])\geq d(x,u) \geq d(x)/2.$$
The proof of this lemma is complete.
\epf

\begin{lem}\label{z2}
If $\gamma$ is an $a$-cone arc in $X$ joining $x$ to $x_0$, then for each $y\in\gamma[x,x_0]$ and $z\in \gamma[y,x_0]$, we have
\be\label{zz1} k(y,z)\leq \ell_k(\gamma[y,z])\leq 3a\log\left (1+\frac{ad(z)}{d(y)}\right )\leq 3ak(y,z)+3a\log(3a).
\ee
\end{lem}

\bpf Since $\gamma$ is an $a$-cone arc, for each $u\in\gamma[y,z]$, we have
$\ell(\gamma[y,u])\leq a d(u).$  Moreover, Lemma \ref{t-1} ensures that $d(y)\leq 2a d(u).$
These two facts guarantee that
$\ell(\gamma[y,u])+d(y)\leq 3a d(u).$
Therefore,
\beq\nonumber
k(y,z)&\leq &\ell_k(\gamma[y,z]) = \int_{\gamma[y,z]} \frac{|du|}{d(u)}
\\ \nonumber &\leq&  \int_0^{\ell(\gamma[y,z])} \frac{3a dt}{t+d(y)}= 3a\log\left (1+\frac{\ell(\gamma[y,z])}{d(y)}\right )
\\ \nonumber &\leq& 3a\log\left (1+\frac{a d(z)}{d(y)}\right )
\\ \nonumber &\leq& 3a\log\frac{d(z)}{d(y)}+3a\log(3a)\leq 3a k(y,z)+3a\log(3a).
\eeq
\epf

\textbf{Proof of Theorem \ref{thm-0}.}
\emph{Necessity}: Suppose that $X$ is quasihyperbolic $a$-John with center $x_0$. Fix $x\in X$ and take a quasihyperbolic geodesic $\gamma$ in $X$ between $x$ and $x_0$. The assumption ensures that $\gamma$ is an $a$-cone arc. Hence $ad(x_0)\geq \ell(\gamma[x,x_0])\geq d(x,x_0)$ and thus $2ad(x_0)\geq \diam (X)$, by the arbitrariness of $x\in X$.

Moreover, suppose that $y_1$ and $y_2$ are points on $\alpha$ satisfying for each $u\in \alpha[y_1,y_2]$, $d(u)\leq 2\min\{d(y_1),d(y_2)\}$. Therefore, we have $d(y_1)\leq 2d(y_2)$. The last inequality combined with Lemma \ref{z2}, shows that
$$k(y_1,y_2)\leq 3a\log\Big(1+\frac{a d(y_1)}{d(y_2)}\Big)\leq 3a\log(1+2a).
$$
The necessity follows by taking $A=3a\log(1+2a)+2a$.

\emph{Sufficiency}: Suppose that there is a constant $A\geq 1$ such that $A\, d(x_0)\geq \diam (X)$ and for every $z\in X$ and  that for each quasihyperbolic geodesic $\alpha$ between $z$ and $x_0$, the following holds: If $y_1$ and $y_2$ are points on $\alpha$ satisfying for each $u\in \alpha[y_1,y_2]$, $d(u)\leq 2\min\{d(y_1),d(y_2)\}$, then we have $k(y_1,y_2)\leq A.$
Without loss of generality, we may assume that $A\geq 2$.
We need to show that every quasihyperbolic geodesic $\alpha$ in $X$ joining $z$ to $x_0$ is an $a$-cone arc for some $a\geq 1$.

Fix $z\in X$ and take a quasihyperbolic geodesic $\alpha$ in $X$ between $z$ and $x_0$. Let $y_0\in \alpha$ be a point satisfying
$d(y_0)=\max\{d(u): \;u\in \alpha\}.$
Then there is a unique nonnegative integer $n$ such that
$2^nd(z)\leq d(y_0)<2^{n+1}d(z).$
For each $i=0,\ldots,n$, let $z_i$ be the first point on $\alpha$ such that
\be\label{0.1} d(z_i)=2^id(z),\ee
when traveling from $z$ to $x_0$. Similarly, we define $x_j$ to be the first point on $\alpha$ such that
$d(x_j)=2^jd(x_0)$
for $j=0,\ldots,m$ when traveling from $x_0$ to $z$, where $m$ is the unique nonnegative integer such that
$2^md(x_0)\leq d(y_0)<2^{m+1}d(x_0).$

Put $z=z_0$. We observe that the arc $\alpha$ has been divided into $(n+m+1)$ non-overlapping (modulo end points) subarcs
$$\alpha[z_0,z_1],\ldots,\alpha[z_{n-1},z_{n}],\;\alpha[z_n,x_m],\alpha[x_m,x_{m-1}], \ldots,\alpha[x_1,x_0].
$$
Moreover, all subarcs in the above are also quasihyperbolic geodesics between their respective endpoints, and for any one of the above, denoted by $\alpha[x,y]$, we have
$d(u)\leq 2\min\{d(x),d(y)\}$ for all $u\in\alpha[x,y].$
Therefore, for each $0\leq i\leq n-1$, we may apply the assumption to the subarc $\alpha[z_{i},z_{i+1}]$ and obtain by (\ref{zt-2}) that
\be\label{1.0}\log \frac{\ell(\alpha[z_i,z_{i+1}])}{d(z_i)}\leq k(z_i,z_{i+1})\leq A.\ee
This inequality implies that
\be\label{1.1} \ell(\alpha[z_i,z_{i+1}])\leq e^A d(z_i).\ee
Similarly, for each $0\leq j\leq m-1$, we get
\be\label{1.2} \;\;\;k(x_j,x_{i+1})\leq A,\;\;\ell(\alpha[x_{j+1},x_j])\leq e^A d(x_j)\;\;\;\mbox{and}\;\;\;\; \ell(\alpha[z_n,x_m])\leq e^A d(z_n).\ee

To show that $\alpha$ is an $a$-cone arc, we need to prove that $\ell(\alpha[z,v])\leq a\,d(v)$ for each $v\in \alpha$. For this, we consider three cases in terms of the location of $v$ on $\alpha$.

\bca Suppose $v\in \alpha[z_i,z_{i+1}]$ for some $0\leq i\leq n-1$.\eca
By (\ref{zt-1}) and (\ref{1.0}), we have
$$\Big|\log \frac{d(v)}{d(z_i)}\Big|\leq k(z_i,v)\leq k(z_i,z_{i+1})\leq A,$$
which ensures
\be\label{1.3} d(z_i)\leq e^A d(v).\ee
This inequality, together with (\ref{0.1}) and (\ref{1.1}), guarantees that
$$\ell(\alpha[z,v])\leq \sum_{s=0}^{i}\ell(\alpha[z_s,z_{s+1}])\leq e^A \sum_{s=0}^{i} d(z_s)\leq 2e^Ad(z_i)\leq 2e^{2A} d(v).
$$

\bca\label{t-2} Suppose $v\in \alpha[z_n,x_m]$.\eca

A similar argument as (\ref{1.3}) gives $d(z_n)\leq e^Ad(v)$. Using this and (\ref{1.2}), we have
\beq\nonumber
\ell(\alpha[z,v]) &\leq & \sum_{s=0}^{n}\ell(\alpha[z_s,z_{s+1}]) +\ell(\alpha[z_n,x_m])
\\ \nonumber &\leq& 2e^{2A}d(z_n)+e^Ad(z_n)
\leq 3e^{2A}d(v).
\eeq

\bca Suppose $v\in \alpha[x_{j+1},x_j]$ for some $0\leq j\leq m-1$.\eca
By Case \ref{t-2}, we have $\ell(\alpha[z,x_m])\leq 3e^{2A}d(x_m).$  Since $Ad(x_0)\geq \diam (X)\geq d(y_0)$, we have $m\leq 2^m\leq A$. Using (\ref{zt-2}) and (\ref{1.2}), we obtain
$$\log\frac{d(x_m)}{d(v)}\leq \log\left (1+\frac{\ell(\alpha[x_0,x_m])}{d(x_m)}\right )\leq k(x_0,x_m)=\sum_{j=0}^{m-1}k(x_j,x_{j+1}) \leq mA\leq A^2
$$
which ensures that
$\ell(\alpha[x_0,x_m])\leq e^{A^2}d(x_m)$ and $d(x_m)\leq e^{A^2}d(v).$
Therefore
$$\ell(\alpha)=\ell(\alpha[z,x_m])+\ell(\alpha[x_0,x_m])\leq (3e^{2A}+e^{A^2})d(x_m)\leq 4e^{2A^2}d(v).
$$
Hence the quasihyperbolic geodesic $\alpha$ is an $a$-cone arc with $a=4e^{2A^2}$, which completes the proof of Theorem \ref{thm-0}.
\qed

\subsection{Gromov hyperbolic John spaces are quasihyperbolic John} In this part, we prove Theorem \ref{thm-1} by using Theorem \ref{thm-0}.  In fact, we use
Lemma \ref{z2} and Theorem \ref{z1} to verify the criterion for quasihyperbolic John spaces as follows.

\begin{lem}\label{z3}
Let $(X,d)$ be an $a$-John space with center $x_0\in X$, and let $(X,k)$ be
Gromov $\delta$-hyperbolic, where $k$ is the quasihyperbolic metric of $X$. Then
there is a constant $A\geq 1$ depending only on $a$ and $\delta$  such that $A d(x_0)\geq \diam (X)$ and for every $x\in X$ and for each quasihyperbolic geodesic $\alpha$ between $x$ and $x_0$, the following holds: If $y_1$ and $y_2$ are points on $\alpha$ satisfying for each $u\in \alpha[y_1,y_2]$, $d(u)\leq 2\min\{d(y_1),d(y_2)\}$, then we have
$k(y_1,y_2)\leq A.$
\end{lem}

\bpf Since $X$ is $a$-John with center $x_0$, for every $x\in X$ there is an $a$-cone arc $\gamma$ connecting $x$ and $x_0$ in $X$. Thus
$ad(x_0)\geq \ell(\gamma)\geq d(x,x_0)$ and so $2ad(x_0)\geq \diam (X)$, by the arbitrariness of $x\in X$.

Next, fix $x\in X$ and pick a quasihyperbolic geodesic $\alpha$ between $x$ and $x_0$. The assumption implies that there is also an $a$-cone arc $\gamma$ connecting $x$ and $x_0$ in $X$. Using Lemma \ref{z2}, $\gamma$ is actually  a quasihyperbolic $(\lambda,\mu)$-quasigeodesic with $\lambda=3a$ and $\mu=3a\log(3a)$. Since $(X,k)$ is
Gromov $\delta$-hyperbolic, Theorem \ref{z1} ensures that there is   an $R=R(\lambda,\mu,\delta)=R(a,\delta)$ such that the quasihyperbolic Hausdorff distance
satisfies the inequality
\be\label{3.1} k_\mathcal{H}(\alpha,\gamma)\leq R. \ee

Let $y_1$ and $y_2$ be the points on $\alpha$ satisfying $d(u)\leq 2\min\{d(y_1),d(y_2)\}$ for all $u\in \alpha[y_1,y_2]$. We need to estimate $k(y_1,y_2)$. Using (\ref{zt-1}) and (\ref{3.1}), we see that there are $z_1$ and $z_2$ on $\gamma$ such that for $i=1,2$,
\be\label{3.2}\Big|\log \frac{d(y_i)}{d(z_i)}\Big|\leq k(y_i,z_i)\leq R,\ee
which implies
$e^{-R}d(y_i)\leq d(z_i)\leq e^{R}d(y_i).$
Since $2^{-1}d(y_1)\leq d(y_2)\leq 2d(y_1)$, we obtain  $2^{-1}e^{-R}d(z_1)\leq d(z_2)\leq 2e^R d(z_1)$. It follows from Lemma \ref{z2} that
$$k(z_1,z_2)\leq 3a\log(1+2ae^R),
$$
 which, combined with (\ref{3.2}), ensures that
$$k(y_1,y_2)\leq 2R+3a\log(1+2ae^R).
$$
Therefore, Lemma \ref{z3} is true by taking $A=2a+2R+3a\log(1+2ae^R)$.
\epf

\textbf{Proof of Theorem \ref{thm-1}.}  Theorem \ref{thm-1} follows from Lemma \ref{z3} and Theorem \ref{thm-0}.
 \qed


\end{document}